\documentclass[secthm,seceqn,amsthm,ussrhead]{amsart}
\usepackage{amsmath,latexsym}
\usepackage[psamsfonts]{amssymb}
\usepackage{times}
\usepackage[mathcal]{euscript}
\numberwithin{equation}{section}

\newcommand{\bbT}{\mathbb T}

\renewcommand{\epsilon}{\varepsilon}

\newcommand{\be}{\begin{equation}}
\newcommand{\ee}{\end{equation}}
\newcommand{\no}{\nonumber}

\newcommand{\C}{\mathbb{C}}

\newcommand{\R}{\mathbb{R}}

\newcommand{\T}{\mathbb{T}}

\newcommand{\Z}{\mathbb{Z}}

{\bf}{\it}
\newtheorem{theorem}{Theorem}[section]
\newtheorem{lemma}[theorem]{Lemma}
\newtheorem{corollary}[theorem]{Corollary}

\newtheorem{definition}[theorem]{Definition}

\newtheorem{remark}[theorem]{Remark}

\date{\today}

\begin{document}
\title{Low energy effects for a family of  Friedrichs models under rank one perturbations}

\author{Sergio  Albeverio$^{1,2,3}$, Saidakhmat  N. Lakaev$^{4,5}$,
  Ramiza Kh. Djumanova $^{5}$}

\address{$^1$ Institut f\"{u}r Angewandte Mathematik,
Universit\"{a}t Bonn, Wegelerstr. 6, D-53115 Bonn\ (Germany)}

\address{
$^2$ \ SFB 611, \ Bonn, \ BiBoS, Bielefeld - Bonn;}
\address{
$^3$ \ CERFIM, Locarno and Acc.ARch,USI (Switzerland) E-mail
albeverio@uni.bonn.de}

\address{
{$^4$ Samarkand State University,Samarkand (Uzbekistan)} \ {E-mail:
lakaev@yahoo.com }}

\address{
$^5$ {Samarkand division of Academy of sciences of
Uzbekistan(Uzbekistan)}}

\maketitle
\begin{abstract}

A family of Friedrichs models with rank one perturbations
$h_\mu(p),$ $p \in (-\pi,\pi]^3,\mu>0$ associated to a system of two
particles on the lattice $\Z^3$ is considered. The existence of a
unique strictly positive eigenvalue below the bottom of the
essential spectrum of $h_\mu(p)$ for all nontrivial values $p \in
(-\pi,\pi]^3$ under the assumption that $h_\mu(0)$ has either a zero
energy resonance (virtual level) or a threshold eigenvalue is
proved. Low energy asymptotic expansion for the Fredholm determinant
associated to family of Friedrichs models is obtained.
\end{abstract}



\section{Introduction}
In the present paper we consider a family of Friedrichs models under
rank one perturbations associated to a system of two particles on
the lattice $\Z^3.$

The main goal of the paper is to give a thorough mathematical
treatment of the spectral properties of a family of Friedrichs
models with emphasis on low energy expansions for the Fredholm
determinants associated to the family (see, e.g.\cite{AGH,ALMM,
ALzMahp,FIC,Lfa93,Sob,Tam94,Yaf2} for relevant discussions and
\cite{K-M,M-S},
 \cite{Zh}
for the general study of the low-lying excitation spectrum for
quantum systems on lattices).

These kind of models have been discussed in quantum mechanics
\cite{Fad,Fri}, solid physics \cite{RSIII, Mat, Mog,G-Sh} and in
lattice field theory \cite{MalMin,LMfa04,LMfa05}.

Threshold energy resonances (virtual levels) for the two-particle
Schr\"odinger operators have been studied in
\cite{AGH,ALzMahp,ALMM,K-S,Lfa93,Yaf2}. The threshold expansions for
the resolvent of two-particle Schr\"odinger operators have been
studied in \cite{ALzMahp,J-K,K-S,Ltmf92,Lfa93,Sob,Tam94,Yaf2} and
have been applied to the proof of the existence of Efimov's effect
in \cite{ALzMahp,Lfa93,Sob,Tam94,Yaf1}.

Similarly to the lattice  Schr\"odinger operators and in contrast to
the continuous Schr\"odinger operators the family of Friedrichs
models $h_\mu(p), p \in (-\pi,\pi]^3,\mu>0$ depends parametrically
on the internal binding $p $, the quasi-momentum, which ranges over
a cell of the dual lattice and hence it has spectral properties
analogous to those of lattice Schr\"odinger operators.

Let us recall that the spectrum and resonances of the original
Friedrichs model and its generalizations have been studied  and the
finiteness of the eigenvalues lying below the bottom of the
essential spectrum has been proven in \cite{Fri,Fad,Ltsp86,Yaf3}.

In \cite{LMfa04,LMfa05} a peculiar family of Friedrichs models was
considered and the appearance of eigenvalues for values of the total
quasi-momentum $p \in (-\pi,\pi]^d,d=1,2$ of the system lying in a
neighbourhood of some particular values of the parameter $p$ has
been proven.

For a wide class of the two-particle  Schr\"odinger operators the
existence of eigenvalues of $h_\mu(p), p \in (-\pi,\pi]^3$ for all
nonzero values of the quasi-momentum $0\neq p\in \T^3$ has been
proven in \cite{ALMM}.

The main results of the present paper are as follows.

First of them gives the existence of a unique strictly positive
eigenvalue $e_\mu(p)$ below the bottom of the essential spectrum of
$h_\mu(p), p \in (-\pi,\pi]^3$ for all nonzero values of the
quasi-momentum $0\neq p\in \T^3$ (provided that $h_\mu(0)$ has
either a zero energy resonance or a threshold eigenvalue) and upper
bound on it (Th. \ref{pos.eig}).

The second one allows monotonous dependence of the eigenvalue
$e_\mu(p)$ on $\mu$ (Th. \ref{mon.eig}).

The third one presents an expansion for the Fredholm determinant
resp. the Birman-Schwinger operator in powers of the quasi-momentum
$p$ in a small $\delta$-neighborhood of the origin and proving that
the Fredholm determinant resp.the Birman-Schwinger operator has an
analytic  continuation to the bottom of the essential spectrum of
$h_{\mu}(p)$ as a function of $w=(m-z)^{1/2}\geq 0$ for $z\leq m,$
where $z\in \R^1$ is spectral the parameter (Th. \ref{main}).

The structure of the paper is as follows. In Sec. 2 states the
problem and presents the main results. Proofs are presented in Sec.
4 and are based on a series of lemmas in Sec. 3.

 Throughout the present paper we adopt the following
conventions:
 $\T^3$ denotes  the three-dimensional
torus, identified with the cube $(-\pi,\pi]^3$ with appropriately
identified sides. For each sufficiently small $\delta>0$ the
notation $U_{\delta}(0) =\{p\in {\bbT}^3:|p|<\delta \}$ stands for a
$\delta$-neighborhood of the origin. Denote by $L_2(\Omega)$ the
Hilbert space of square-integrable functions defined on a measurable
subset $\Omega$ of $\R^n.$

\section{The family of operators  $h_\mu(p),\,  p\in {\T}^3$ and
statement of main results}

We introduce the following family of bounded self-adjoint operators
(the Friedrichs model) $h_\mu(p),\,  p\in {\T}^3$ acting in
$L_2(\T^3)$ by
\begin{equation}\label{h_alpha}
 h_\mu(p)=h_{0}(p)-\mu v.
 \end{equation}
The non-perturbed operator $h_{0}(p)$ is the multiplication operator
by the function $u(p,q)$
\begin{equation}
 (h_0(p)f)(q)=u(p,q)f(q),\quad
 f\in L_2(\T^3),\,\\
\end{equation}
\begin{equation}
u(p,q)\equiv\varepsilon (p)+ \varepsilon (p-q)+ \varepsilon (q),
\end{equation}
\begin{equation}
\varepsilon (p)\equiv 3-cos p_1-cos p_2 -cos p_3,\,p=(p_1,p_2,p_3)
\in {\T}^3
\end{equation}
\begin{remark}
The operator $h_0(p)$ is associated to a system of two particles
(bosons) moving on the three-dimensional lattice $\Z^3$ and is
called the {\it\bf  free Hamiltonian}.
\end{remark}

The perturbation operator is an integral operator of rank one and
has the form
\begin{align*}
 \label{poten} (vf)(q)=\varphi(q)\int\limits_{\T^3}
\varphi(t)f(t)dt,
 f\in L_2(\T^3),\\
\end{align*}
where $\varphi$ is real-analytic even function defined on ${\T}^3$
and $\mu>0$ is a positive real number.

Since the perturbation $v$ of the multiplication operator $h_0(p)$
is a self-adjoint integral operator of rank one in accordance to
Weyl's theorem on invariance of the essential spectrum under finite
rank perturbations the essential spectrum of the operator $h_\mu(p)$
fills  the following interval on the real axis:
$$
\sigma_{ess}(h_\mu(p))=[m(p),M(p)],\,\,
$$
where
$$
m(p)=\min_{q\in {\T}^3}u(p,q) ,\,\,M(p)=\max_{q\in {\T}^3} u(p,q).
$$
\begin{remark}\label{minimum}
Notice that the function $u(\cdot,\cdot)$ defined on $({\T}^3)^2$
can be represented in the form
 \begin{equation}\label{}
u(p,q)=\varepsilon(p)+\sum_{i=1}^3(2-2cos\frac{p_i}{2}cos(q_i-\frac{p_i}{2})),\quad
p,q\in (-\pi,\pi)^3.
\end{equation}
and hence for any $p { \in } (-\pi,\pi)^3$  the point
  $q_0(p)=p/2\in (-\pi,\pi)^3$ is the unique non-degenerate minimum
of the function $u_p(\cdot)=u(p,\cdot).$
\end{remark}
\begin{remark}
For some $p\in \T^3$ (for example $p=(\pi,\pi,\pi)\in \T^3$) the
essential spectrum of $h_\mu(p)$ can  degenerate to the set
consisting  of the unique point $[m(p), M(p)].$ Because of this we
can not state that the essential spectrum of $h_\mu(p)$ is
absolutely continuous for any $p\in \T^3.$
\end{remark}

Let $C(\T^3)$ be the Banach space of continuous functions on
$\T^3.$

\begin{definition}\label{resonance0}
 The operator $h_\mu(0)$ is said to
have a virtual level at the bottom of the essential spectrum (zero
energy resonance) if the number $1$ is an eigenvalue of the
integral operator
$$
(\mathrm{G_\mu}\psi)(q)=\mu \varphi(q) \int\limits_{{\T}^3}
(u(0,t))^{-1}\varphi(t)\psi(t)dt,\,\,\psi\in {C(\T^3)}$$ and the
associated eigenfunction $\psi$ satisfies the condition $\psi(0)\neq
0.$
\end{definition}
\begin{remark} We note that the operator $h_\mu(0)$ has a zero energy resonance
if and only if $\varphi(0)\neq 0$ and $\mu=\mu_0$(see Lemma
\ref{resonance}), where
$$
 \mu_0=\Big( \int\limits_{{\T}^3}
\varphi^2(t) (u(0,t))^{-1}dt \Big)^{-1}.
$$

\end{remark}
\begin{remark}
(i) If $\varphi(0)\not=0$ and the operator $h_{\mu}(0)$ has a zero
energy resonance, then the function
\begin{equation}\label{f0f1}
f(q)=\varphi(q)u_0(q))^{-1},
\end{equation}
obeys the equation $ h_{\mu}(0)f=0$ and  $ f\in L_1(\T^3)\setminus
L_2(\T^3)$
(see Lemma \ref{resonance}).\\
(ii) If $\varphi(0)=0$ and zero is an eigenvalue of the operator
$h_{\mu}(0),$ then the function $f$ defined by \eqref{f0f1}, obeys
the equation $ h_{\mu}(0)f=0$ and $ f\in L_2(\T^3)$ (see Lemma
\ref{zeroeigen}).
\end{remark}

The following theorems will be proven in Sect. 4 and are based on
lemmas proven in Sect. 3.
\begin{theorem}\label{pos.eig}
For all nonzero $p\in \T^3$ the operator $h_{\mu_0}(p),\,p\in \T^3$
has a unique strictly positive eigenvalue $e_{\mu_0}(p).$ One has
$$0<e_{\mu_{0}}(p)<m(p),\,0\neq p\in \T^3.$$
\end{theorem}
\begin{theorem}\label{mon.eig}
 For any $\mu>\mu_0$ the operator $h_\mu(p),\,p\in \T^3$ has a
unique eigenvalue $e_\mu(p).$ One has
$$e_\mu(p)<e_{\mu_0}(p)<m(p),\,0\neq p\in
\T^3$$ and
$$e_\mu(0)<0.$$
\end{theorem}

Let ${\C}$ be the field of complex numbers. For any $p \in \T^3$ we
define an analytic  function $\Delta_{\mu}(p,\cdot)$(the Fredholm
determinant
 associated to the operator $h_\mu(p)$)  in
 ${\C} { \setminus } [m(p),M(p)]$ by
\begin{equation*}\label{det}
\Delta_{\mu}(p,z)=1-\mu
\int\limits_{{\T}^3}(u(p,t)-z)^{-1}\varphi^2(t)dt.
\end{equation*}

\begin{theorem}\label{main}
(i) For any  $z\leq m(p)$ the function $\Delta_\mu(\cdot,z)$ is
analytic on $U_{\delta}(0)$  and the following decomposition
\begin{equation}
\Delta_\mu(p,z)=\Delta_\mu(0,z)+\Delta_\mu^{res}(p,z)\,,p\in
U_\delta(0)
\end{equation}
holds, where $\Delta_\mu^{res}(p,z)=O(|p|^2)$ as $p\to 0$ uniformly
in $z\leq m(p).$

(ii)  For  $\Delta_\mu(0,\cdot)$  the following decomposition
\begin{equation}\label{D(0zeta)-1}
\Delta_\mu(0,z))=\Delta_\mu(0,0)- \pi^2 \mu
\varphi^2(0)(-z)^{-1/2}+\Delta_\mu^{res}((z))
\end{equation}
holds, where $\Delta_\mu^{res}(z))=O(|z|)$ as $z\to -0$ and
$(-z)^{-1/2}>0,$ for $z<0$
\end{theorem}

Theorem \ref{main} and Lemma \ref{maximum} imply the  following
corollary on the asymptotic behavior of $\Delta_{\mu}(p,0)$ near
$p=0,$ which plays an important role in the proof of the finiteness
or infiniteness of the number of eigenvalues (Efimov's effect) for a
model operator associated to a system of three-particles on the
lattice $\Z^3$ \cite{ALzMahp}.

\begin{corollary}\label{razlojeniya} (i) Let the operator $h_{\mu_0}(0)$ have a
zero energy resonance. Then for some $c_1,c_2>0$  the following
asymptotics
\begin{align}\label{raz}
&\Delta_{\mu_0}(p,0)= \frac{\sqrt{3}}{2}\pi^2 \mu_0
\varphi^2(0)|p|+O(|p|^2)\quad\mbox{as}\quad p\to 0
\end{align}
and inequalities  hold
\begin{equation}c_1
|p|\le \Delta_{\mu_0}(p,0) \le c_2 |p|,\, p\in U_{\delta}(0).
\end{equation}

(ii)Let the operator $h_{\mu_0}(0)$ have a zero eigenvalue. Then for
some $c>0$  the following inequality  holds
\begin{align*}
\Delta_{\mu_0}(p,m)\geq c p^2,
p\in U_\delta(0).\\
\end{align*}

\end{corollary}

\section{Some spectral properties of the operators  $h_\mu(p),\,  p\in {\T}^3$ }
In this section we study some spectral properties of the family
$h_\mu(p),\,  p\in {\T}^3$ with emphasis on the zero energy
resonance and zero eigenvalue.
\begin{lemma}\label{delta=0}
For any $\mu>0$ and $p\in \T^3$ the following statements are
equivalent:\\ (i) The operator $h_\mu(p)$ has an eigenvalue $z
\in {\C} \setminus [m(p),M(p)]$ below the bottom of the essential spectrum.\\
(ii) $\Delta_{\mu}(p,z)=0$, $z \in
{\C} \setminus [m(p),M(p)].$\\
 (iii) $\Delta_{\mu}(p,z')<0$ for some
$z'\leq m(p).$
\end{lemma}
\begin{proof}
From the positivity of $v$ it follows that there exists a unique
positive square  root of $v,$ which will be denoted
$v^{\frac{1}{2}}.$ For any $\mu>0$ and $p\in \T^3$ the number $z \in
{\C} \setminus [m(p),M(p)]$ is an eigenvalue of $h_\mu(p)$ if and
only if  the number $\lambda=1$ is an eigenvalue of the operator
\begin{equation*}
G_{\mu}(p,z)=\mu v^{\frac{1}{2}} r_0(p,z)v^{\frac{1}{2}}
\end{equation*}
(this follows from the Birman-Schwinger principle). Since the
operator $v^{\frac{1}{2}} $ is of the form $$ (v^{\frac{1}{2}}
f)(q)=||\varphi||^{-1}\varphi(q)\int\limits_{\T^3}
\varphi(t)f(t)dt,\quad f \in L_2(\T^3) $$ the operator
$G_{\mu}(p,z)$ has the form
\begin{equation*}
(G_{\mu}(p,z)f)(q)=\mu ||\varphi||^{-2}{\Lambda}(p,z)
\varphi(q)\int\limits_{\T^3} \varphi(t)f(t)dt,\quad f \in L_2(\T^3),
\end{equation*}
where
\begin{equation}\label{Lambda}
{\Lambda}(p,z)=\int\limits_{{\T}^3} (u(p,t)-z)^{-1}\varphi^2(t)dt.
\end{equation}

According to  the Fredholm theorem  the number $\lambda=1$ is an
eigenvalue for the operator $G_{\mu}(p,z)$ if and only if
$\Delta_{\mu}(p,z)=0. $ The equivalence of $(i)$ and $(ii)$ is
proven.

Now we prove the equivalence of $(ii)$ and $(iii)$. Let
$\Delta_{\mu}(p,z_0)=0$ for some $z_0 \in \C\setminus
[m(p),M(p)]$. The operator $h_\mu(p)$ is self-adjoint and hence
the equivalence $(i)$ and $(ii)$ implies that $z_0$ is a real
number. From $\Delta_{\mu}(p,z)>1$ for all $z>M(p)$ we conclude
that $z_0\in (-\infty,m(p)).$
 Since the function $\Delta_{\mu}(p,\cdot),\,p\in
{\T}^3$ is decreasing in $z\in (-\infty,m(p))$ we have
$\Delta_{\mu}(p,z')< \Delta_{\mu}(p,z_0)=0$ for some $
z_0<z'<m(p).$

Now we suppose that $\Delta_{\mu}(p,z')<0$ for some $z' \leq m(p).$
For any $p\in {\T}^3$ the function $\Delta_{\mu}(p,\cdot)$ is
continuous in $z\in (-\infty,m(p))$ and  $\lim\limits_{z\to -\infty}
\Delta_{\mu}(p,z)=1,$  hence  there exists $z_0\in (-\infty;z')$
such that $\Delta_{\mu}(p,z_0)=0.$ This completes the proof.
\end{proof}

The function $u(p,\cdot),p\in (-\pi,\pi)^3$ has a unique
non-degenerate minimum at $q=p/2$ (see Remark \ref{minimum}) and
hence by Lebesgue's dominated convergence theorem  the finite limit
$$
\Delta_{\mu}(p,m(p))=\lim_{z\to m(p)-} \Delta_{\mu}(p,z),p\in
{(-\pi,\pi)}^3.
$$ exists.

Note that
$$
\Delta_{\mu}(\pi,z)=1-\mu
(12-z)^{-1}\int\limits_{{\T}^3}\varphi^2(t)dt
$$
and
$$\lim_{z\to 12-} \Delta_{\mu}(\pi,z)=+\infty.$$
The following lemmas  describe whether the bottom of the essential
spectrum of $h_{\mu_0}(0)$ is a zero energy resonance or eigenvalue.
\begin{lemma}\label{resonance}   The following statements are equivalent:\\
(i) The operator $h_{\mu}(0)$ has a zero energy resonance.\\
(ii)  $\varphi(0)\neq 0$  and $\Delta_{\mu}(0, 0)=0.$\\
(iii) $\varphi(0)\neq 0$  and $\mu= \mu_{0}.$
\end{lemma}

\begin{proof} We
prove\,\,$(i)\rightarrow(ii)\rightarrow(iii)\rightarrow(i).$ Let the
operator $h(0)$ have a zero energy  resonance for some $\mu>0$. Then
by
 Definition \ref{resonance0} the equation
\begin{equation}\label{res-def}
\psi(q)=\mu \varphi(q) \int\limits_{{\T}^3}
(u(0,t))^{-1}\varphi(t)\psi(t)dt,\,\,\psi\in {C(\T^3)}
 \end{equation}
 has a nontrivial solution $\psi\in C({\T^{3}}),$ which satisfies the condition
 $\psi(0)\neq 0.$

In fact,  this solution is equal to the  function $\varphi$ (up to
a constant factor) and hence $\Delta_{\mu}(0,0)=0$
 and so
$\mu=\mu_{0}$.

 Let for some $\mu>0$ the equality
$\mu=\mu_{0}$ hold and consequently $\Delta_{\mu}(0,0)=0$. Then the
function $\varphi\in  C({\T}^{3})$  (up to a constant factor) is a
solution of the equation \eqref{res-def},
  that is, the
operator $h_{\mu_0}(0)$ has a zero energy resonance.
\end{proof}

\begin{lemma}\label{zeroeigen}
  The following statements are equivalent:\\
(i) the operator $h_\mu(0)$ has a zero eigenvalue.\\
(ii)  $\varphi(0)=0$  and $\Delta_{\mu}(0,0)=0.$\\
(iii) $\varphi(0)=0$ and $\mu= \mu^{0}.$
\end{lemma}
\begin{proof} We
prove\,\,$(i)\rightarrow(ii)\rightarrow(iii)\rightarrow(i).$ Suppose
that $f\in L_2(\T^3)$ is an eigenfunction of the operator $h(0)$
associated to a zero eigenvalue. Then $f$ satisfies the equation
\begin{equation}\label{h=0}
u(0,q)f(q)-\mu \varphi(q) \int\limits_{{\T}^3} \varphi(t)f(t)dt=0.
\end{equation}
From \eqref{h=0} we find that $f,$ except for an arbitrary factor,
is given by
\begin{equation}\label{eigenunction}
f(q)=(u(0,q))^{-1}\varphi^2(q),
\end{equation}
and from \eqref{h=0} we derive the equality $\Delta_{\mu}(0,0)=0.$

The functions $u(\cdot,\cdot)$ and $\varphi(\cdot)$ are analytic on
$(\T^3)^2$ and $\T^3,$ respectively and the function $u(0,\cdot)$
has a unique non-degenerate minimum at the point $q=0$ and hence
$f\in L_2(\T^3)$ implies that $\varphi(0)=0.$

Substituting the expression \eqref{eigenunction} for $f$ to the
equation \eqref{h=0} we get the equality
\begin{equation*}
\varphi(q)=\mu \varphi(q) \int\limits_{{\T}^3}
(u(0,t))^{-1}\varphi^2(t)dt.
 \end{equation*}

Hence $\Delta_{\mu}(0,0)=0$ and so $\mu=\mu_0.$

Let $\varphi_{\alpha}(0)=0$ and  $\mu=\mu_0,$ then
 $\Delta_{\mu}(0,0)=0,$ and the function
$f,$ defined by \eqref{eigenunction}, obeys the equation $h(0)f=0$
and $ f\in L_2(\T^3).$
 \end{proof}
\begin{lemma}\label{maximum}
The function $\Lambda(\cdot,\cdot),$ defined by \eqref{Lambda}, has
a non-degenerate maximum at the point $p=0$.
\end{lemma}
\begin{proof}  Since
$u(\cdot,\cdot)$ and $\varphi(\cdot)$ are even functions on
$(\T^3)^2$ and $\T^3,$ respectively, the function $\Lambda(\cdot,0)$
is also even on $\T^3.$ Then we get

\begin{align}\label{L.Even}
\Lambda(p,0)-\Lambda(0,0)\no\\
&=\frac{1}{4}\int\limits_{{\T}^3} \frac{2u(0,t)-(u(p,t)+u(-p,t))}
{u(p,t)u(-p,t)u(0,t)}[u(p,t)+u(-p,t)]\varphi^2(t)dt\no\\
&-\frac{1}{4}\int\limits_{{\T}^3} \frac{ [u(p,t)-u(-p,t)]^2}
{u(p,t)u(-p,t)u(0,t)}\varphi^2(t)dt.
\end{align}
For all $\mid p\mid+\mid t\mid\neq 0$ we have $u(p,t)+u(-p,t)>0.$
Moreover the inequality
\begin{equation*}
u(0,t)-\frac{u(p,t)+u(p,-t)}{2}= \sum_{i=1}^{3}(\cos p_i-1)(1+\cos
t_i)>0
\end{equation*}
holds for all $p\neq0$ and $t\neq(\pi,\pi,\pi).$ Hence
\eqref{L.Even} implies the inequality
$$\Lambda(p,0)-\Lambda(0,0)<0$$ for all  $0\neq p\in \T^3,$ that is,
the function $\Lambda(\cdot,0)$ has a unique maximum at $p=0.$

The equalities
\begin{align*}
&  \frac{\partial } {\partial p^{(i)} }u_0(t) =\sin t^{(i)},\,\,
\frac{\partial^2 } {\partial p^{(i)}\partial p^{(i)}
}u_0(t) =1+\cos t^{(i)},\\
& \frac{\partial^2 } {\partial p^{(i)}\partial p^{(j)}
}u_0(t)=0,\quad i\neq j,\,\,i,j=1,2,3,\nonumber
\end{align*}
yield the  relations
\begin{align}\label{par.der0}
&\frac{\partial^2 \Lambda(0,0)} {\partial p^{(i)}
\partial p^{(i)}}<0,
\frac{\partial^2 \Lambda(0,0)} {\partial p^{(i)}
\partial p^{(j)}}=0,\,\,i\neq j,\,\,i,j=1,2,3
\end{align}
and hence using \eqref{par.der0} we get that the matrix of the
second order partial derivatives of the function $ \Lambda(\cdot,0)
$ at the point $p=0$ is negative definite. Thus the function
$\Lambda(\cdot,0)$ has a non-degenerate maximum at the point $p=0.$
\end{proof}

 Set
$$
\C_+=\{z\in \C: Re\, z>0\},\quad \R_+=\{x\in \R: x>0\},\quad
\R_+^0=\R_+\cup \{0\}.
$$
Let $u_0(\cdot,\cdot)$ be the function defined on
$U_{\delta}(0)\times \T^3\,$ by
 \begin{equation}\label{w}
u_0(p,q)=u_p(q+p/2))-m(p),
\end{equation}
where $q=p/2,\,p\in (\pi,\pi)^3$ is the non-degenerate minimum point
of the function $u_p(\cdot)$ (see Remark 2.2.).

For any $p\in U_\delta(0)$ we define an analytic function
$D(p,\cdot)$ in $\C_+$ by
\begin{equation*}
D(p,w)= \int\limits_{{\T}^3}
(u_0(p,q)+w^2)^{-1}{\varphi}^2(q+p/2)dq.
\end{equation*}

\begin{lemma}\label{raz}
 (i) For any  $w\in \C_+$ the function
$D(\cdot,w)$ is analytic in $U_{\delta}(0)$  and the following
decomposition
\begin{equation}\label{D(pzeta)-1}
D(p,w)=D(0,w)+D^{res}(p,w)\,,p\in U_\delta(0)
\end{equation}
holds, where $D^{res}(p,w)=O(p^2)$ as $p\to 0,$ uniformly in $w\in
\R_+^0.$

(ii)  The derivative of $D(0,\cdot)$ at $w=0$ exists and the
decomposition
\begin{equation}\label{D(0zeta)-1}
D(0,w)=D(0,0)+ \pi^2 \mu \varphi^2(0)w+D^{res}(w)
\end{equation}
holds, where $D^{res}(w)$ is analytic function and hence
$D^{res}(w)=O(w^{2})$ as $w\to 0.$
\end{lemma}
\begin{proof}
 Since $m(\cdot)$
is analytic in $U_{\delta}(0)$ by definition of the function
$D(\cdot,\cdot)$ we obtain that the function $D(\cdot,w)$ is
analytic in $U_{\delta}(0)$ for any $w\in \C_+.$

The asymptotics
 \begin{align*}
u_0(p,q)=q^2+O(|p|^2|q|^2)+O(|q|^4) \,\,as\,\,|p|,|q| \to 0
\end{align*}
yields  the existence of $C>0$ such that for any $w\in R^{0}_{+}$
and $i,j=1,2,3$ the following inequalities hold
\begin{align}\label{Estimate1}
&\Big |\frac{\partial^2}{\partial p_i \partial p_j}
 \frac{\varphi(q+q_0(p))}{u_0(p,q)+w^2}\Big |\leq \frac{C}{q^2},\,p,q\in U_\delta(0)
\end{align}
and
\begin{align}\label{Estimate2}
&\Big |\frac{\partial^2}{\partial p_i \partial p_j}
 \frac{\varphi(q+q_0(p))}{u_0(p,q)+w^2}\Big |\leq C,\,p\in U_\delta(0), q\in \T^3
 \setminus U_\delta(0).
\end{align}
 Then Lebesgue's dominated convergence theorem implies
$$
\frac{\partial^{2}}{\partial p_i\partial p_j}D(p,0)=\lim_{w\to
0+}\frac{\partial^{2}}{\partial p_i\partial p_j}D(p,w).
$$

Repeatedly applying Hadamard's lemma (see \cite{Zor} V.1, p. 512) we
obtain
\begin{align*}
D(p,w)=D(0,w)+\sum_{i=1}^3\left[\frac{\partial}{\partial p_i}
D(p,w)\right]_{p=0} p_i+\sum_{i,j=1}^{3}H_{ij}(p,w) p_i p_j,
\end{align*}
where  the functions $H_{ij}(\cdot,w),\,w\in \R_+^0,\,\,i,j=1,2,3$
are continuous in $U_\delta(0)$ and
\begin{equation*}
    H_{ij}(p,w)=\frac{1}{2}\int_{0}^1\int_{0}^1
    \frac{\partial^2}{\partial p_i \partial p_j}D(x_1x_2p,w)dx_1dx_2.
\end{equation*}

The estimates \eqref{Estimate1} and \eqref{Estimate2} give
$$
|H_{i,j}(p,w)|\leq \frac{1}{2}\int_{0}^1\int_{0}^1
    \Big |\frac{\partial^2}{\partial
p_i \partial p_j}D(x_1x_2p,w)\Big |dx_1dx_2\leq
C\Big(1+\int\limits_{U_\delta(0)}\frac{dq}{q^2}\Big)
$$
 uniformly in  $p\in U_\delta(0)$ and $w\in \R_+^0.$

Since for any $w\in \R_+^0$ the function $D(\cdot,w)$ is even in
$U_\delta(0)$ we have
$$
\left[\frac{\partial}{\partial p_i} D(p,w)\right]_{p=0}
=0,\,i=1,2,3.
$$

${\bf(ii)}$ The function  $D(0,\cdot)$ can be analytically continued
to $\C_+\cup V_\gamma(0)$ (see \cite{Ltmf92}),where $V_\gamma(0)$ is
a ball of the radius $\gamma>0$ with the center at $w=0$ in $\C.$
Denote by $D^*(0,\cdot)$ this analytic continuation. Then  the
representation
\begin{equation*}
D^*(0,w)=D(0,0)+ \frac{\partial}{\partial w}D(0,0)w+D^{res}(w)
\end{equation*}
holds, where  $D^{*,res}(w)$ is analytic function in $V_\gamma(0)$
and hence $D^{*,res}(w)=O(w^{2})$ as $ w\to 0.$

Now we prove the  equality
\begin{equation}\label{part D0}
\frac{\partial}{\partial w}D(0,0)=\pi^2 \mu \varphi^2(0).
\end{equation}
The function  $u_0(0,\cdot)$ has a unique non-degenerate minimum at
$q=0.$ Therefore,  by virtue of the Morse lemma (see \cite{Zor})
there exists a one-to-one mapping $\psi:W_\gamma(0)\rightarrow
\tilde W(0)$ of a certain ball $W_\gamma(0)\subset \T^3$  of radius
$\gamma>0$ with the center at $t=0$ to a neighborhood $\tilde
W(0)\subset\T^3$ of the point \,$q=0$ such that:
\begin{align}\label{eps=t2}
u_0(0,\psi(t))=t^2
\end{align}
with $\psi(0)=0$ and the equality $J_\psi(0)=1/2$ holds, where
$J_\psi(t)$  is the Jacobian of the mapping $q=\psi(t).$ For any
$w\in \C_+$ the function $\frac{\partial}{\partial w}D(0,\cdot)$ can
be represented in the form
\begin{align}\label{partialD}
\frac{\partial}{\partial w}D(0,w)=D^{(1)}(w)+D^{(2)}(w),\quad w\in
\C_+
\end{align}
with
\begin{align}\label{D=1}
D^{(1)}(w)=2w\int\limits_{\T^3\setminus \tilde W(0)}
(u_0(0,q)+w^2)^{-1}\varphi^2(q)dq ,\quad w\in \C_+
\end{align}
and
\begin{align}\label{D=2}
D^{(2)}(w)=2w\int\limits_{\tilde W(0)}
(u_0(0,q)+w^2)^{-1}\varphi^2(q)dq ,\quad w\in \C_+.
\end{align}

Since the function $u_0(0,\cdot)$ is continuous on the compact set
$\T^3\setminus \tilde W(0)\subset\T^3$ and has a unique minimum at
$q=0$ there exists $c>0$ such that $|u_0(0,q)|>c$ for all $q\in
\T^3\setminus U_{\delta}(0).$

Then we have
\begin{align}\int\limits_{\T^3\setminus \tilde W(0)}
(u_0(0,q)+w^2)^{-1}\varphi^2(q)dq\rightarrow
\int\limits_{\T^3\setminus \tilde
W(0)}(u_0(0,q)+0)^{-1}\varphi^2(q)dq
\end{align}
 as $w \rightarrow 0.$ In the integral \eqref{D=2} making a change
of variable $q=\psi(t)$ and using the equality \eqref{eps=t2} we
obtain
\begin{align}\label{D1-D1}
D^{(2)}(w)=2w\int\limits_{W_\gamma(0)}
(t^2+w^2)^{-1}\varphi^2(\psi(t))J_\psi(t)dt.
\end{align}

Going over in the integral in \eqref{D1-D1} to spherical coordinates
$t=r\omega,$ we reduce it to the form
\begin{equation}\label{I2}
D^{(2)}(w)=2w \int_0^{\gamma}r^2(r^2+w^2)^{-1}F(r) dr
\end{equation}
with
$$
F(r)=\int_{\Omega_2}\varphi^2(\psi(r\omega))J_\psi(r\omega)d\omega,
$$
where $\Omega_2$ is the unit sphere in $\R^3$ and $d \omega$ is the
element of the unit sphere in this space.

It easy to see that
\begin{equation*}\label{Limit}
\int_{0}^{\gamma}
  {w}(r^2+w^2)^{-1}dr \to \frac{\pi}{2},\quad \mbox{as}\ w \to 0+
\end{equation*}
and
\begin{equation*}
\int_{0}^{\gamma}
  w r^2(r^2+w^2)^{-1} (F(r)-F(0))dr\to 0,
  \quad \mbox{as}\ w \to 0+.
\end{equation*}
Hence the equality \eqref{part D0} holds.
\end{proof}

\section{The proof of the main results}

{\it\bf{The proof of Theorem \ref{pos.eig}}}. It easy to see that
for all $0\neq p\in \T^3$ and a.e. $ q \in \T^3$ the inequality
$u_p(q)-m(p)<u_0(q)$ holds. Hence for all $0\neq p\in \T^3$ the
inequality
$$\Lambda(p,m(p))-\Lambda(0,0)>0$$
yields
$$
\Delta_{\mu_0}(p,m(p))<\Delta_{\mu_0}(0,0)=0,\,p\in \T^3.$$ By Lemma
3.4 we have
$$\Delta_{\mu_0}(p,0)>\Delta_{\mu_0}(0,0)=0,\,p\in
\T^3.$$ Since $ \lim_{z\to -\infty} \Delta_{\mu_0}(p,z)=1$ and
$\Delta_{\mu_0}(p,\cdot)$ is monotonously decreasing
 on $(-\infty,m(p))$ we conclude that the
function $\Delta_{\mu_0}(p,z)$ has a unique solution in $(0,m(p)).$
Lemma \ref{delta=0} completes the proof of Theorem \ref{pos.eig}.

{\it\bf{Proof of Theorem \ref{mon.eig}}}. Let $\mu>\mu_0.$ We have
$$ \Delta_{\mu}(p,z)<\Delta_{\mu_0}(p,z)$$ for all $ p\in \T^3,z\leq m(p).$
By Theorem 2.7 for all nonzero $p\in \T^3$ the operator
$h_{\mu_0}(p),\,p\in \T^3$ has a unique strictly positive eigenvalue
$e_{\mu_0}(p)$ and it satisfies the inequality
$$0<e_{\mu_{0}}(p)<m(p),\,0\neq p\in
\T^3.$$ Lemma \ref{delta=0} yields
$\Delta_{\mu_0}(p,e_{\mu_0}(p))=0,\,0\neq p\in \T^3.$ Since the
function $\Delta_{\mu}(p,\cdot)$ is decreasing on $(-\infty,m(p))$
 for all $
0\neq p\in \T^3,z< m(p)$ we get
$$\Delta_{\mu}(p,z)< \Delta_{\mu_0}(p,e_{\mu_0}(p))= 0$$  and
$$\Delta_{\mu}(0,z)< \Delta_{\mu_0}(0,0)= 0,z<0.$$
Taking into account $ \lim_{z\to -\infty} \Delta_{\mu}(p,z)=1$ and
applying Lemma \ref{delta=0} again we complete the proof of Theorem
\ref{pos.eig}.$\Box$

{\it\bf{The proof of Theorem \ref{main}}} follows from Lemma
\ref{raz} if we take into account that $w=(m(p)-z)^{1/2}\geq 0$ for
$z\leq (m(p).$$\Box$

{\bf Acknowledgement} This work was supported by the DFG 436 USB
113/3 project and the Fundamental Science Foundation of Uzbekistan.
The second named author gratefully acknowledge the hospitality of
the Institute of Applied Mathematics and of the IZKS of the
University of Bonn.

\end{document}